\documentclass{article}
\usepackage[T1]{fontenc}
\usepackage{graphicx, a4wide}
\usepackage{amssymb,amsfonts,amsmath,amsthm,mathtools} 
\usepackage{subfigure} 
\usepackage{ae,aecompl}     

\newcommand\unnumberedfootnote[1]{ %
        \let\temp=\thefootnote %
        \renewcommand{\thefootnote}{}%
        \footnote{#1}%
        \let\thefootnote=\temp%
        \addtocounter{footnote}{-1}}

\theoremstyle{definition}

\numberwithin{equation}{section}

\makeatletter
\renewcommand{\@fnsymbol}[1]{\@arabic{#1 }}
\makeatother

\begin{document}

\title{Asymmetric autocatalytic reactions and their stationary distribution}

\author{Cameron Gallinger, Lea Popovic\thanks{Department of
    Mathematics and Statistics Concordia University, Montreal QC H3G
    1M8, Canada; email: lea.popovic@concordia.ca}}

\date{November 27. 2023}

\maketitle

\unnumberedfootnote{\emph{Keywords and phrases.} reaction networks, autocatalytic reactions, stationary distribution, discreteness induced transitions, Moran model, genic selection}

\begin{abstract}  We consider a general class of autocatalytic reactions, that has been shown to display stochastically switching behaviour (Discreteness Induced Transitions) in some parameter regimes. This behaviour was shown to occur when either the overall species count is low, or when the rate of inflow and outflow of species is relatively much smaller than the rate of autocatalytic reactions.
The long-term behaviour of this class was analyzed in \cite{KBW20} with an analytic formula for the stationary distribution in the symmetric case. We focus on the case of asymmetric autocatalytic reactions and provide a formula for an approximate stationary distribution of the model. We show this distribution has different properties corresponding to the distinct behaviour of the process in the three parameter regimes. In the Discreteness Induced Transitions regime, the formula provides the fraction of time spent at each of the stable points. 
 \end{abstract}

\section{Introduction} Many important processes in biology rely on the switching behaviour between distinct states of the internal state of the system. Mathematically, this may be the result of alternating between different stable states for the dynamics of the molecular composition of the system. This feature is present in many gene-expression systems, where a gene alternates between two types of states (“on” and “off”) regulating the production of a protein. It is also present in many phosphorylation switches in signalling pathways. Such bistable switching patterns arise from distinct stable equilibria in the deterministic dynamics and the ability of infrequent large stochastic fluctuations to pull the system from a basin of attraction of one equilibrium to the other. However, there are also cases of chemical dynamics in which bistability is not possible in the deterministic model and is only possible in the stochastic model of the same chemical reaction system (e.g. \cite{BQ10}). The importance of stochastic effects in biological switching is well recognized and a variety of stochastic processes are being used to model them (see \cite{B17} for a review).

Autocatalytic reaction networks are an important class of reaction networks whose long-term behaviour can substantively differ when modelled deterministically versus stochastically (e.g. \cite{VQ07}). Autocatalytic reactions are of broad interest in modelling living systems as they involve only a small number of distinct species and interactions, but they exhibit a range of outcomes such as self-sustaining growth, oscillations and symmetry breaking. Autocatalysis occurs in many elements of cellular metabolism including glycolysis, mitosis, apoptosis, and DNA replication (for a review see \cite{HSNS21}), and may plausibly play a key role in the origin of life (e.g.\cite{HS18}). They were recently used in chemical computing,  \cite{AC21},  as part of an artificial neural network design to perform digital image recognition tasks.

We analyze a model for a simple example of an autocatalytic reaction network   proposed by Togashi-Kaneko in \cite{TK01}. They were particularly interested in the dynamics of the network when the total species count in the system was low, and observed numerically that stochasticity at the low species count led to a new type of switching dynamics. The system spent most of its time in states where some combination of species is absent, switching rapidly by a successive increase in some species count from one such state to another. They called this behaviour Discretely Induced Transitions (DITs) because discrete changes in species counts drive the switching. 
This dynamics is prominent in the parameter regime when the rates of autocatalytic reactions are fast relative to the rates of inflow and outflow from a reservoir, and the average time between switches is  large. 

The long-term behaviour of this model was further explored  in \cite{TK03}, where the same authors observed that, if only the long term averages are observed, autocatalytic reaction networks with strong symmetries do not show the important distinction between small and large species counts. However, in systems possessing asymmetry, numerical simulations clearly show that the shape of the long-term distributions shifts as the species count decreases. They observed similar effects in simple catalytic reaction networks in  \cite{TK07}, where the low count of some species type effectively acts by switching on and off parts of the reaction network. They also observed this effect in spatial models in  \cite{TK05}, where higher abundances are typically modelled by reaction-diffusion systems, but the localization and low abundance of some species types generates different possibilities for the long term behaviour patterns of the system. 
 
First analytic results for symmetric auto-catalytic reaction networks were obtained by a stochastic differential equation approximation of the system dynamics in \cite{BRM12}. They used a scaling parameter composed of the product of the system volume and inflow rate, and derived an analytic expression for the stationary distribution of a linear combination of the species in the system, based on a time-scale separation. In \cite{BDM14} they used the same tools to explore a version of this auto-catalytic reaction network without inflows and outflows, but with single species conversions. They provided a formula for the mean switching time between the states when the system contains only one of the species. A rigorous approximation by an obliquely reflected stochastic differential equation of a general autocatalytic reaction network was subsequently developed in \cite{FYY22}, where a constrained Langevin approximation (based on the theory of  Leite and Williams \cite{LW19}) was analyzed,  together with the stationary distribution of this reflected stochastic differential equation.

 Stochastic approximations are sensitive to the different abundances in molecular amounts (e.g. \cite{BKPR06}), as well as to the time-scale separations in the rates of different reactions in the network (e.g. \cite{McSP14}). For this reason, it is worth analyzing the properties of reaction networks from the Markov chain model directly.
 Analytic results for the Markov chain model were explored in \cite{SK15}, where the effect of species counts was shown on the flow of the molecular 'current'  from one species type to another. Analytic results for the stationary distribution of a related class of models were derived in \cite{HM19}, for reaction networks without inflows and outflows and hence conserved overall species counts. 
 A thorough analytic exploration of the auto-catalytic reaction network with inflows and outflows was performed in \cite{KBW20}, where it was shown that the Markov chain model of a general auto-catalytic reaction network is exponentially ergodic. The authors further analyzed the symmetric version of this model and derived an explicit form for the stationary distribution as a Dirichlet-multinomial distribution (Theorem 4.3, \cite{KBW20}). However, an explicit formula for the case of asymmetric autocatalytic reaction networks remained unknown.  

To analyze the stationary distribution in the asymmetric case, we explored the connection of the auto-catalytic model with inflows and outflows to that of the Moran model with mutation from population genetics. A connection with the Moran model was noted in  \cite{BDM14, SK15} for symmetric rates of auto-catalysis. We observe that asymmetric rates appear in Moran models that contain genic selection, as in \cite{EG09}, and the formula for the stationary distribution for such a Moran model is the starting point in our analysis. We explore the consequences due to the differences between the two models and explore how close the stationary distribution for the Moran model is to that of the auto-catalytic network. We analyze the effect of the asymmetry in the long-term distribution of the species count in the system. We show the approximate stationary distribution also has the signature of the Discretely Induced Transitions in the parameter regime corresponding to low species counts, or slow inflow and outflow rates  in the system.




\section{Stochastic (Markov chain) model for auto-catalytic reactions}


The auto-catalytic reaction model we examine is the Togashi-Kaneko model for an open system. It consists of $d$ molecular species $\{A_i\}_{ i=1,\dots, d}$ involved in auto-catalytic reactions and subject to inflows and outflows from the system:
\begin{equation}\label{acreactions}
    A_i +A_j \xrightarrow{\kappa_{i}} 2A_i \quad\quad \emptyset\xrightleftharpoons[\delta]{\lambda_i} A_i\qquad i,j=1,2,\dots,d
\end{equation}
The stochastic model is a continuous-time Markov chain $(X(t))_{t\ge 0}$, where $X(t)=(X_1(t),\!\dots\!, X_d(t))$ counts the number of molecules of species $i\in\{1,\dots, d\}$ present in the system at time $t\ge 0$. The state space for the Markov chain is $E=\{a=(a_1,\dots,a_d)\in \mathbb{N}^d\}$ and the transition rates are given by mass-action kinetics as
\begin{equation}\label{transrates}
    \begin{split}
     &q_{a,a-e_j+e_{i}}=\kappa_{i} a_ia_j,\qquad \forall i\neq j\in\{1,\dots, d\}\\
      &q_{a,a+e_i}=\lambda \quad q_{a,a-e_i}=\delta a_i, \quad \forall i\in\{1,\dots, d\}
    \end{split}
\end{equation}
Note that this is not necessarily a cyclic model (as all pairs of different species can catalyze each other), and the rates of autocatalytic reactions are different for each pair. 
The only constraint on the reactions that we impose is that the outflow rate is the same for all species, which ensures the lumpability of the Markov chain into a process that counts the total mass in the system. 

The concept of {\it lumpability} means that the state space $E$ can be partitioned into subsets $E_n=\{a\in\mathbb{N}^d|\sum_{i=1}^d=n\}$, and that the {\it lumped} process $(N(t))_{t\ge 0}$ defined by $N(t)=\sum_{i=1}^d X_i(t)$ is a continuous time Markov chain  on $E_n$ with associated transition rates
\begin{equation}
    q_{n,n+1}=\sum_{i=1}^d \lambda_i \quad\quad q_{n,n-1}=n\delta.
\end{equation}
Because the total mass process $(N(t))_{t\ge 0}$ is a linear death process with immigration on $\mathbb{Z}_+$, we know that it is irreducible and positive recurrent with stationary distribution $\nu$ that is Poisson with mean equal to $\sum_{i=1}^d\lambda_i/\delta$. 

In \cite{KBW20} Theorem 4.1,  Bibbona, Kim and Wiuf showed that the Markov chain $(X(t))_{t\ge 0}$ corresponding to \eqref{transrates} is positive recurrent on $E$ with a unique stationary distribution $\Pi$, and that it converges exponentially fast to this distribution. 

In the special case when all the auto-catalytic rates are the same: $\kappa_i=\kappa, \forall i$, they derive an explicit formula for the stationary distribution $\Pi$ in the form 
\begin{equation}\label{Pi}
\Pi(a)=\pi(a|n)\nu(n)
\end{equation}
where $\nu\sim$ Poisson$(\mu)$ is the stationary distribution of the total mass process, and $ \pi(a|n)$ is the Dirichlet-multinomial$(n,\alpha)$ distribution on $E_n$  
\begin{equation}\label{bibbonapi}
   \pi(a|n)={n\choose a} \frac{\Gamma\big(\sum_{i=1}^d\alpha_i\big)}{\Gamma\big(n+\sum_{i=1}^d\alpha_i\big)}\prod_{i=1}^d \frac{\Gamma\big(a_i+\alpha_i\big)}{\Gamma\big(\alpha_i\big)}
\end{equation}
with parameters $\alpha=(\alpha_1,\dots, \alpha_d)$ given by
\begin{equation}\label{alpha}
\alpha_i=\frac{\delta\lambda_i}{\kappa \sum_{j=1}^d\lambda_j}\end{equation}
This analytic formula for $\Pi$ from \eqref{Pi} and \eqref{bibbonapi} is directly verified using the time evolution for the probability distribution of the system.

\subsection{Connection to population genetics}
The Moran model is a fundamental stochastic model for the evolution of a population of species types undergoing reproduction and mutation in continuous time (overlapping generations). The model describes a population of fixed size $n$, with  $d$ different  species types $\{A_i\}_{ i=1,\dots, d}$ , whose number of individuals at time $t\ge 0$ is denoted by $\tilde X(t)=(\tilde X_1(t),\dots, \tilde X_d(t))$.  
Each individual of type $i$ gives birth at the rate $\kappa_i$, and its offspring replaces an individual chosen at random to die. The species undergo mutations with each individual of type $i$ changing to type $j$ at rate of $v p_{ij}$, where $p_{ij}$ are the entries of a stochastic matrix $P$ and $v>0$ is the overall mutation rate. 

The species reproduction and mutation changes can be represented graphically by: 
\begin{equation}
    A_i+A_j \xrightarrow{\kappa_j} 2A_j \quad\quad A_i \xrightarrow{v p_{j}}A_j
\end{equation}
where we assumed parent-independent mutations, $p_{ij}=p_j $. 
The fact that the fitness rates $\kappa_j$ are not equal across species is due to genic selection, and the special case  $\kappa_i=\kappa, \forall i\in\{1,\dots,d\}$ corresponds to a model that is neutral under selection.

 The continuous time Markov chain of species counts $(\tilde X(t))_{t\ge 0}$ has  transition rates from $a$ to $a-e_i+e_j$ equal to
\begin{equation}
    \tilde q_{a,a-e_i+e_j}=a_i\bigg[\frac{\kappa_j}{n}a_j+\nu p_{j}\bigg].
\end{equation}
The Moran model differs from the model of auto-catalytic reactions \eqref{acreactions} in that it has a conservation of mass (closed system) and lacks the (open system) effects of the inflow and outflow reactions. 

A comparison of long-term behaviour for closed versus open systems exists for deterministic models of reaction networks: if the open system is created from the closed system by adding inflow and outflow for each species type, then the open system has the same number of steady states and stability as the closed system (\cite{CF06, BP18, CFW20}). For stochastic systems, the comparison is different as the closed system may have absorbing states when the open system does not affect the existence and the number of modes in the stationary distribution. The reason the Moran model is remarkably close to the TK model in its long-term behaviour is due to the fact that one can create a correspondence between the mutation events in the former and the in-/out-flows in the latter (see explanation above \eqref{alphaTK}).


\section{Approximate stationary distribution}
The Moran model with genic selection on a population of size $n$ was shown by Etheridge and Griffiths \cite{EG09} to have a stationary distribution proportional to 
\begin{equation}
    \pi_n(a)\propto \kappa_1^{a_1}\cdots\kappa_d^{a_d}{n\choose a} \frac{\alpha_{1_{(a_1)}}\cdots\alpha_{d_{(a_d)}}}{|\alpha|_{(n)}}, \; a\in E_n
\end{equation}
where the parameters are given by
\begin{equation}\label{alphaMoran}
    \alpha_i=\frac{nv p_i}{\kappa_i},
\end{equation}
$|\alpha|=\sum_{i=1}^d\alpha_i$, and $\alpha_{(a)}$ is the ascending factorial notation (Pochammer function)
\[\alpha_{(a)} =\alpha(\alpha+1)···(\alpha+a-1).\] The normalizing constant for $\pi_n(a)$ is given as the partition function
\begin{equation}\label{partitionmoran}
    u(\alpha,\kappa,n)=\mathbb{E}\bigg[\bigg(\sum_{i=1}^d\kappa_i\xi_i\bigg)^n\bigg]
\end{equation}
where $\xi=(\xi_1,\xi_2,...,\xi_d)$ has a Dirichlet-multinomial$(n,\alpha)$ distribution. Writing the ascending factorial in terms of Gamma functions $\alpha_{(a)}=\Gamma(\alpha+a)/\Gamma(\alpha)$, the stationary distribution can also be written as 
\begin{equation}\label{moranpi}
    \pi_n(a)\propto
     \prod_{i=1}^d \kappa_i^{a_i} {n\choose a} \frac{\Gamma(\sum_{i=1}^d\alpha_i)}{\Gamma(\sum_{i=1}^d \alpha_i +n)}\prod_{i=1}^d \frac{\Gamma(\alpha_i+a_i)}{\Gamma(\alpha_i)}
\end{equation}

In case of the neutral model, when $\kappa_i=\kappa, \forall i=1,\dots,d$, the partition function becomes $u(\alpha,\kappa,n)=\kappa^n\mathbb{E}\left[\left(\sum_{i=1}^d \xi_i\right)^n\right]=\kappa^n$ since $\sum_{i=1}^d \xi_i=1$. The stationary distribution then reduces to the Dirichlet-multinomial distribution
\begin{equation*}\label{multidir}
   \pi_n(a)={n\choose a} \frac{\Gamma(\sum_{i=1}^d\alpha_i)}{\Gamma(\sum_{i=1}^d \alpha_i +n)}\prod_{i=1}^d \frac{\Gamma(\alpha_i+a_i)}{\Gamma(\alpha_i)}
\end{equation*}
which is the same as the stationary distribution $\pi(a|n)$ in \eqref{bibbonapi} derived by Bibbona, Kim and Wiuf  under the parameter change \[\alpha_i=\frac{\delta \lambda_i}{\kappa \sum_{j=1}^d \lambda_j}.\] 

To create a correspondence between the Moran model with genic selection and the auto-catalytic reactions model, suppose we require that each outflow of some species is identified in time with the next inflow of some species (e.g. $A_i$), thus closing the system. With this modification the overall rate of inflow of species $A_i$ (simultaneously with some outflow) would be $\delta \lambda_i/(\sum_{j=1}^d \lambda_j)$. In the Moran model the overall rate of mutation into species $A_i$ is $nvp_i$. Identifying these inflow rates in the equation for the parameter \eqref{alphaMoran} would result in using the parameter 
 \begin{equation}\label{alphaTK}
 \alpha_i=\frac{\delta \lambda_i}{\kappa_i \sum_{j=1}^d \lambda_j}
 \end{equation}
 for the asymmetric auto-catalytic reaction system.

The correspondence suggests the distribution $\pi_n$ from \eqref{moranpi} may be a good candidate in forming the approximate stationary distribution for asymmetric autocatalytic reaction system, and we use $\tilde\pi_n$  to denote $\pi_n$ with parameters $\alpha_i$ satisfying \eqref{alphaTK}. Our proposed approximate stationary distribution for the autocatalytic reaction system \eqref{acreactions} then is 
\begin{equation}\label{tildepi}
    \tilde\Pi(a)=\nu(n)\frac{1}{u(\alpha,\kappa,n)}
     \bigg[\prod_{i=1}^d \kappa_i^{a_i}\bigg] {n\choose a} \frac{\Gamma(\sum_{i=1}^d\alpha_i)}{\Gamma(\sum_{i=1}^d \alpha_i +n)}\prod_{i=1}^d \frac{\Gamma(\alpha_i+a_i)}{\Gamma(\alpha_i)}, \;\; n=\sum_{i=1}^da_i
\end{equation}
We next analyze how well the distribution $\tilde \Pi$ approximates  the true stationary distribution $\Pi$. 

\subsection{Analytical results}

To assess how close $\tilde \Pi$ is to the true stationary distribution of \eqref{acreactions}, let $\mathcal{A}$ denote the generator of the Markov chain  $(X(t))_{t\ge 0}$ defined for any  $f:E\mapsto \mathbb{R}$ by
\begin{align}\label{generator}
{\mathcal A}f(a)&=\sum_{i,j=1}^d q_{a,a-e_j+e_i}[f(a-e_j+e_i)-f(a)] \\&+\sum_{i=1}^dq_{a,a+e_i}[f(a+e_i)-f(a)]+\sum_{i=1}^dq_{a,a-e_i}[f(a-e_i)-f(a)]\nonumber
\end{align}
where the transition rates $q_{i,j}$ are as in \eqref{transrates}. Let  $p_t(\mu,a)$ denote the probability that the process $(X(t))_{t\ge 0}$ is in state $a\in E$ at time $t\ge 0$, given that its initial distribution is $X(0)\sim \mu$. Let ${\mathcal A}^*$ be the adjoint of ${\mathcal A}$ defined 
by  
\begin{align}\label{adjoint}
{\mathcal A}^*p_t(\mu,a)&=\sum_{i\neq j=1}^d \big[q_{a+e_j-e_i,a}\, p_t(\mu,a+e_j-e_i)\mathbf{1}_{a+e_j-e_i\in E} -  q_{a,a+e_j-e_i}p_t(\mu,a)\big]\\
&+\sum_{i=1}^d\big[q_{a+e_i,a}p_t(\mu,a+e_i)-q_{a,a+e_i}p_t(\mu,a)\big]+\sum_{i=1}^d\big[q_{a-e_i,a}p_t(\mu,a-e_i)-q_{a,a-e_i}p_t(\mu,a)\big]\nonumber
\end{align}
which specifies the evolution of $p_t(\mu,a)$ via the Kolmogorov forward equation ({\it chemical master equation} in the biology literature) as 
\begin{equation}\label{Kolmogorov}
\frac{d}{dt} p_t(\mu,a) = {\mathcal A}^* p_t(\mu,a)\mathbf{1}_{a\in E}.
\end{equation}
The stationary distribution $ \Pi$ of \eqref{acreactions} is in the null space of ${\mathcal A}^* $, satisfying the {\it global balance equation}: ${\mathcal A}^* \Pi(a)=0$.  Using the proposed distribution $\tilde\Pi=\nu(n)\tilde\pi_n(a)$ equation \eqref{adjoint} becomes
\begin{align}\label{adjointpi}
{\mathcal A}^*\tilde\Pi(a)&=\nu(n)\sum_{i\neq j=1}^d \big[q_{a+e_j-e_i,a}\, \tilde\pi_n(a+e_j-e_i)-  q_{a,a+e_j-e_i}\tilde\pi_n(a)\big]\\
&+\nu(n)\sum_{i=1}^d\big[q_{a+e_i,a}\frac{\sum_i\lambda_i}{\delta(n+1)}\tilde\pi_{n+1}(a+e_i)-q_{a,a+e_i}\tilde\pi_n(a)\big]\nonumber\\
&+\nu(n)\sum_{i=1}^d\big[q_{a-e_i,a}\frac{\delta \,n}{(\sum_i\lambda_i)}\tilde\pi_{n-1}(a-e_i)-q_{a,a-e_i}\tilde\pi_n(a)\big]\nonumber
\end{align}
where  $\nu\sim$Poisson$(\sum_i\lambda_i/\delta)$ was used to express $\frac{\nu(n+1)}{\nu(n)}$ and $\frac{\nu(n)}{\nu(n-1)}$ in the second and third sums.

We now use ${\mathcal B}^*(a):={\mathcal A}^*\tilde\Pi(a)/\tilde\Pi(a)$ to analyze how close $\tilde\Pi$ is to $\Pi$. By \eqref{adjointpi} we have
\begin{align*}
{\mathcal B}^*(a)&=\sum_{i\neq j=1}^d \big[q_{a+e_j-e_i,a}\, \frac{\tilde\pi_n(a+e_j-e_i)}{\tilde\pi_n(a)}-  q_{a,a+e_j-e_i}\big]\\
&+\sum_{i=1}^d\big[q_{a+e_i,a}\frac{\sum_i\lambda_i}{\delta(n+1)}\frac{\tilde\pi_{n+1}(a+e_i)}{\tilde\pi_n(a)}-q_{a,a+e_i}\big]
+\sum_{i=1}^d\big[q_{a-e_i,a}\frac{\delta \,n}{(\sum_i\lambda_i)}\frac{\tilde\pi_{n-1}(a-e_i)}{\tilde\pi_n(a)}-q_{a,a-e_i}\big]\nonumber
\end{align*}
For $d=2$ each sum has two terms which we gather (based on $\pm$ sign) into
\begin{align*}
&{\mathcal B}^*(a)=\frac{1}{\tilde\pi_n(a)}\big[ (L_{n-1}+L_n+L_{n+1})(a)-R_n(a)\big]\\
\mbox{with}&\\
    &R_n(a)=[\lambda_1+\lambda_2 +n\delta +(\kappa_1+\kappa_2)a_1a_2]\tilde\pi_n(a)\\
    &L_{n-1}(a)=\frac{n\delta \lambda_1}{\lambda_1+\lambda_2}\tilde\pi_{n-1}(a-e_1) + \frac{n\delta\lambda_2}{\lambda_1+\lambda_2}\tilde\pi_{n-1}(a-e_2)\\
    &L_n(a)=\kappa_2(a_1+1)(a_2-1)\tilde\pi_n(a+e_1-e_2)+\kappa_1(a_1-1)(a_2+1)\tilde\pi_n(a-e_1+e_2)\\
    &L_{n+1}(a)=\frac{\lambda_1+\lambda_2}{n+1}(a_1+1)\tilde\pi_{n+1}(a+e_1)+\frac{\lambda_1+\lambda_2}{n+1}(a_2+1)\tilde\pi_{n+1}(a+e_2)
\end{align*}
We see here the effect of the inflow and outflow reactions that give rise to transitions from $n$ to $n\pm 1$, in different values of $n$ in $\tilde\pi$ and the Poisson distribution $\nu$.
In the expression for $\tilde \pi$, the weighted Dirichlet-multinomial depends only on $a\in E_n$, but the normalizing factor $u$ depends on $n$, so we use Gauss hypergeometric functions to simplify the ratios of $\tilde\pi$ appearing in ${\mathcal B}^*$.

In terms of a Dirichlet-multinomial$(n,\alpha)$ variable $\xi=(\xi_1,\dots,\xi_d)$, the factor $u$ is, by \eqref{partitionmoran}
\begin{equation}
     u(\alpha,\kappa,n)=\mathbb{E}\bigg[\bigg(\sum_{i=1}^d\kappa_i\xi_i\bigg)^n\bigg]
\end{equation}
which we can also write, using 
the product moments of $\xi$, as
\begin{equation*}   
     u(\alpha,\kappa,n)
=\sum_{|a|=n}\bigg[\prod_{i=1}^d\kappa_i^{a_i}\bigg]{n\choose a}\mathbb{E}\bigg[ \prod_{i=1}^d \xi_i^{a_i}\bigg]=
\sum_{|a|=n}\bigg[\prod_{i=1}^d\kappa_i^{a_i}\bigg] {n\choose a}\frac{\Gamma(\sum_{i=1}^d\alpha_i)}{\Gamma(\sum_{i=1}^d\alpha_i+n)}\prod_{i=1}^d \frac{\Gamma(\alpha_i+a_i)}{\Gamma(\alpha_i)}.
\end{equation*}
For $d=2$ this can be further simplified using special functions. Using $a_1=i, a_2=n-i$, 
\begin{equation*}
     u(\alpha,\kappa,n)=\sum_{i=0}^n{n\choose i}\kappa_1^{i}\kappa_2^{n-i}\frac{\Gamma(\alpha_1+\alpha_2)}{\Gamma(\alpha_1+\alpha_2+n)} \frac{\Gamma(\alpha_1+i)}{\Gamma(\alpha_1)}\frac{\Gamma(\alpha_2+n-i)}{\Gamma(\alpha_2)}
\end{equation*}
and the partition function becomes 
\begin{align*}
       u(\alpha,\kappa,n)
&=\frac{\Gamma(\alpha_1+\alpha_2)}{\Gamma(n+\alpha_1+\alpha_2)}\frac{\Gamma(\alpha_2+n)}{\Gamma(\alpha_2)}\kappa_2^n\sum_{i=0}^n{n\choose i} (-1)^i\frac{(\alpha_1)_{(i)}}{(1-\alpha_2-n)_{(i)}}\bigg(\frac{\kappa_1}{\kappa_2}\bigg)^i\\
      &=\frac{\Gamma(\alpha_1+\alpha_2)}{\Gamma(n+\alpha_1+\alpha_2)}\frac{\Gamma(\alpha_2+n)}{\Gamma(\alpha_2)}\kappa_2^n{ }_2F_1\big(-n,\alpha_1;1-\alpha_2-n;\frac{\kappa_1}{\kappa_2}\big)
\end{align*}
where $_2F_1$ is the Gauss hypergeometric function, which when evaluated on a nonpositive integer in the first coordinate reduces to the polynomial
$_2F_1(-n,x;y;z)=\sum_{i=0}^n (-1)^i {n \choose i} \frac{(x)_{(i)}}{(y)_{(i)}}{z^i}$.

\noindent The distribution $\tilde\pi_n$ on $i\in \{0,1,\dots,n\}$ becomes
\begin{align}\label{pid=2}
   \tilde \pi_n(i)&=\frac{1}{u(\alpha,\kappa,n)}
     \kappa_1^{i}\kappa_2^{n-i} {n\choose i} \frac{\Gamma(\alpha_1+\alpha_2)}{\Gamma(\alpha_1+\alpha_2+n)} \frac{\Gamma(\alpha_1+i)\Gamma(\alpha_2+n-i)}{\Gamma(\alpha_1)\Gamma(\alpha_2)}\nonumber\\
    &=\frac{1}{{ }_2F_1\big(-n,\alpha_1;1-\alpha_2-n;\frac{\kappa_1}{\kappa_2}\big)}\bigg(\frac{\kappa_1}{\kappa_2}\bigg)^a {n\choose i}  \frac{\Gamma(\alpha_1+i)\Gamma(\alpha_2+n-i)}{\Gamma(\alpha_1)\Gamma(\alpha_2+n)}.
\end{align}
Our representation by Gauss hypergeometric functions is valid when: (i) $y$ is not a non-positive integer, and when (ii) $ |z| < 1$, or $|z| \le 1$ when $y > x-n$. This requires that: (i) $\alpha_2\notin \mathbb{N}-n$, and that (ii) we label the species so that ${\kappa_1}/{\kappa_2}<1$, allowing for $\kappa_1/\kappa_2\le 1$ when $\delta<\kappa_1$. Since $\alpha_2=(\delta \lambda_2)/\kappa_2(\lambda_1+\lambda_2)$ then $\kappa_1/\kappa_2<1$ together with $\delta<\kappa_1$ will ensure that $0<\alpha_2<1$ as well.\\

Using  \eqref{pid=2} in the expressions $R_n,L_n,L_{n-1},L_{n+1}$ for ${\mathcal B}^*$, we then get that ({\bf see Appendix A}) 
 \begin{align}\label{DBeqn}
     {\mathcal B}^*(a)=&(\lambda_1+\lambda_2)\Bigg[1-\frac{_2F_1(-n,\alpha_1,1-\alpha_2-n,\frac{\kappa_1}{\kappa_2})}{_2F_1(-1-n,\alpha_1,-\alpha_2-n,\frac{\kappa_1}{\kappa_2})}\Bigg]\nonumber\\
     &-\kappa_2(n-1+\alpha_2)\bigg(\frac{a_1\alpha_1}{a_1-1+\alpha_1}+\frac{a_2\alpha_2}{a_2-1+\alpha_2}\bigg)\Bigg[\frac{_2F_1(-n,\alpha_1,1-\alpha_2-n,\frac{\kappa_1}{\kappa_2})}{_2F_1(1-n,\alpha_1,2-\alpha_2-n,\frac{\kappa_1}{\kappa_2})}-1\Bigg]\nonumber\\
      &-\kappa_2\bigg(\frac{a_1\alpha_1}{a_1-1+\alpha_1}+\frac{a_2\alpha_2}{a_2-1+\alpha_2}\bigg)\Big[a_1\big(1-\frac{\kappa_1}{\kappa_2})-\alpha_2\Big]-\kappa_2\frac{a_1\alpha_1}{a_1-1+\alpha_1}\Big[\frac{\kappa_1}{\kappa_2}-1\Big]\nonumber\\
      &-\frac{\lambda_1+\lambda_2}{n+\alpha_2}\frac{_2F_1(-n,\alpha_1,1-\alpha_2-n,\frac{\kappa_1}{\kappa_2})}{_2F_1(-1-n,\alpha_1,-\alpha_2-n,\frac{\kappa_1}{\kappa_2})}\Big[a_1\big(\frac{\kappa_1}{\kappa_2}-1\big)+\alpha_2\Big] 
 \end{align}

We use numerical properties of Gauss hypergeometric functions $_2F_1$ to determine the assumptions on parameters $\lambda_i,\kappa_i,$ and $\delta$ that ensure ${\mathcal B}^*(a)\approx 0$. We also investigate its dependence on the overall species count in the system.
Note that in the symmetric case, $\kappa_1=\kappa_2$, simplifying $R_n,L_{n-1},L_n,L_{n+1}$ yields ${\mathcal B}^*(a)=0$, $\forall a$, and hence $\tilde\Pi$ is the true stationary distribution for $X$ (see proof of Thm 4.2 in \cite{KBW20}). 

\subsection{Scaling and numerical results}

We let $V$ denote the scaling parameter for the magnitude of the overall species count in the system (referred to as {\it volume} in previous analyses of the Togashi-Kaneko model).
The classical mass-action scaling of the stochastic model for the auto-catalytic reaction system \eqref{acreactions} requires all bimolecular reaction rate parameters to scale as $V^{-1}$, all inflow rate parameters to scale as $V$, and all unimolecular outflow rate parameters to stay unscaled,
\[
    \kappa_i=\frac{\kappa_i'}{V},\;\; \lambda_i=\lambda_i'V ,\;\; \delta=\delta'\;\; \quad\Rightarrow \quad \alpha_i=V\alpha_i'
\]

To simplify the expression for ${\mathcal B}^*={\mathcal A}^*\tilde\Pi/\tilde\Pi$, and analyze its properties for low versus high molecular count,
we use the parametrization as in  \cite{KBW20}, with $D,\kappa_i'\sim O(1)$: 
\begin{equation}\label{scalingKBW}
\kappa_i=\frac{\kappa_i'}{V},\;\; \lambda_i=DV,\;\; \delta=D,\;\;  \quad\Rightarrow \quad \alpha_i=\frac{DV}{d\kappa_i'}=V\alpha_i'
\end{equation}
and focus on the role of asymmetric auto-catalytic rates $\kappa_i$ as well as the effect of $V$.
When $d=2$ we rewrite ${\mathcal B}^*$ from \eqref{DBeqn} in terms of scaled parameters ({\bf see Appendix B}) by factoring out the parameters $DV$ and $\kappa_1'/\kappa_2'$ to get
\begin{equation}\label{BV}
 \begin{split}
     {\mathcal B}^*_V(a)=
     &dDV\Bigg(\Bigg[1-\frac{F(-n,\frac{DV}{d\kappa_1'},1-\frac{DV}{d\kappa_2'}-n,\frac{\kappa_1'}{\kappa_2'})}{F(-1-n,\frac{DV}{d\kappa_1'},-\frac{DV}{d\kappa_2'}-n,\frac{\kappa_1'}{\kappa_2'})}\Bigg]\\
     +&\frac{1}{n+\frac{DV}{d\kappa_2'}}\frac{F(-n,\frac{DV}{d\kappa_1'},1-\frac{DV}{d\kappa_2'}-n,\frac{\kappa_1'}{\kappa_2'})}{F(-1-n,\frac{DV}{d\kappa_1'},-\frac{DV}{d\kappa_2'}-n,\frac{\kappa_1'}{\kappa_2'})}\Big[a_1\big(1-\frac{\kappa_1'}{\kappa_2'}\big)-\frac{DV}{d\kappa_2'}\Big]\Bigg)\\
     -&\kappa_2'(n-1+\frac{DV}{d\kappa_2'})\Bigg[\frac{F(-n,\frac{DV}{d\kappa_1'},1-\frac{DV}{d\kappa_2'}-n,\frac{\kappa_1'}{\kappa_2'})}{F(1-n,\frac{DV}{d\kappa_1'},2-\frac{DV}{d\kappa_2'}-n,\frac{\kappa_1'}{\kappa_2'})}-1\Bigg]\bigg(\frac{a_1\frac{DV}{d\kappa_1'}}{a_1-1+\frac{DV}{d\kappa_1'}}+\frac{a_2\frac{DV}{d\kappa_2'}}{a_2-1+\frac{DV}{d\kappa_2'}}\bigg)\\
      -&\kappa_2'\Big[a_1\big(1-\frac{\kappa_1'}{\kappa_2'})-\frac{DV}{d}\Big]\bigg(\frac{a_1\frac{DV}{d\kappa_1'}}{a_1-1+\frac{DV}{d\kappa_1'}}+\frac{a_2\frac{DV}{d\kappa_2'}}{a_2-1+\frac{DV}{d\kappa_2'}}\bigg)+\kappa_2'\frac{a_1\frac{D}{d\kappa_1'}}{a_1-1+\frac{DV}{d\kappa_1'}}\Big[1-\frac{\kappa_1'}{\kappa_2'}\Big].
\end{split}
 \end{equation}
 
First, note that in case $\lambda_i'=\delta'=D\ll \kappa_i'$ then ${\mathcal B}^*_V(a)=O(D)$ implies our proposed distribution $\tilde\Pi$ is close to the stationary distribution, regardless of whether the overall species count $V$ is low or high. Intuitively, when inflow and outflow rates are small relative to the auto-catalytic rates, the system spends most of its time with a constant overall species count, and hence the distribution $\tilde\pi_n$ is stationary most of the time.

\begin{figure}[hbt!]
     \centering
     \includegraphics[width=0.5\textwidth, height=0.195\textheight]{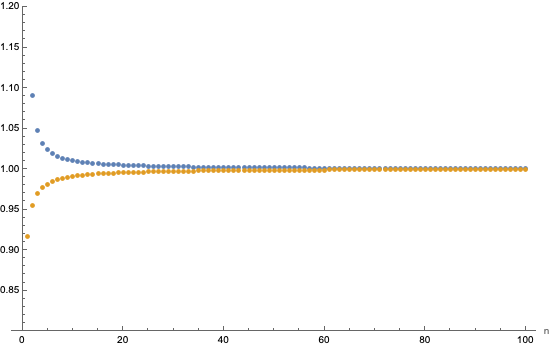}
     \caption{The ratios $\frac{_2F_1(-n,\alpha_1,1-\alpha_2-n,\frac{{\kappa_1}}{{\kappa_2}})}{_2F_1(1-n,\alpha_1,2-\alpha_2-n,\frac{{\kappa_1}}{{\kappa_2}})}$ and $\frac{_2F_1(-n,\alpha_1,1-\alpha_2-n,\frac{{\kappa_1}}{{\kappa_2}})}{_2F_1(-1-n,\alpha_1,-\alpha_2-n,\frac{{\kappa_1}}{{\kappa_2}})}$ plotted in blue and orange, respectively, as a function of n. The rate parameters are $\lambda_1=\lambda_2=2,\,\delta=0.01,\, \kappa_1=1,\, \kappa_2=1.001$}
     \label{fig:hypergeo}
 \end{figure}

\begin{figure}[hbt!]
\begin{subfigure}{}
  \centering
  \includegraphics[width=0.5\linewidth]{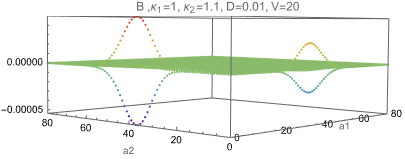}  
  \label{fig:sub-first}
\end{subfigure}
\begin{subfigure}{}
  \centering
  \includegraphics[width=0.495\linewidth]{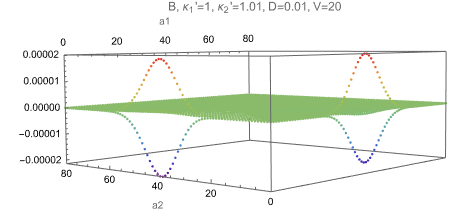}  
  \label{fig:sub-first}
\end{subfigure}
\begin{subfigure}{}
  \centering
  \includegraphics[width=0.495\linewidth]{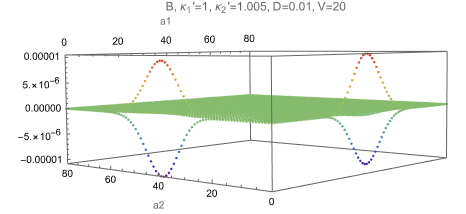}  
  \label{fig:sub-second}
\end{subfigure}
\begin{subfigure}{}
  \centering
  \includegraphics[width=0.5\linewidth]{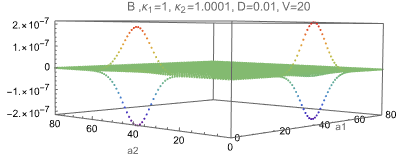}  
  \label{fig:sub-second}
\end{subfigure}
\caption{$\mathcal{B}^*$ plotted with $\lambda_1'=\lambda_2'=\delta'=0.01, V=20$, $\kappa_1'=1$ and $\kappa_2'=1.1, 1.01, 1.005, 1.0001$ gradually decreasing $\kappa_1'/\kappa_2'\to 1$ from top left to bottom right.}
\label{fig:BalanceMultiD}
\end{figure}

Next, numerical explorations show the two ratios of the Gauss hypergeometric functions appearing in ${\mathcal B}^*_V$  are both close to 1 for all $(a_1,a_2)$  (Figure \ref{fig:hypergeo}), so the two terms in the large square brackets are small. Further,
numerical evaluations of ${\mathcal B}^*_V(a)$ for different values of $(a_1,a_2)$ are all close to zero (Figure 2). The very small deviations from zero  occur only when:  $(a_1,a_2) \in \{(0, dV ), (1, dV ), (dV, 0), (dV, 1)\}$. These are values with $a_1+a_2\approx dV$ close to the mean of the overall species count ($\nu\sim$Poisson($\sum_i\lambda_i'/\delta'=dV)$). 
Changing $\kappa'_1/\kappa'_2$ while keeping all other parameters fixed, the deviations from zero disappear as $\kappa'_1/\kappa'_2\to 1$ (see Figure \ref{fig:BalanceMultiD}) as expected.

The effect of the scaling on $\tilde\Pi_V$ is reflected in both $\nu_V\sim$Poisson($\sum_i\lambda_i'/\delta'=dV$) and in the scaling of Gamma functions in  $\tilde\pi_{n,V}$  (the $V$ scaling of $\kappa_i'$ in pre-factor cancels with same in $u_V$) 
\begin{align*}
 &\tilde\pi_{n,V}(a)=\frac{1}{u_V(\alpha',\kappa',n)}{n\choose a}\Bigg[\prod_{i=1}^d(\kappa_i')^{a_i}\Bigg] \frac{\Gamma\big(V\sum_{i=1}^d\alpha_i'\big)}{\Gamma\big(n+V\sum_{i=1}^d\alpha_i'\big)}\prod_{i=1}^d \frac{\Gamma\big(a_i+V\alpha_i'\big)}{\Gamma\big(V\alpha_i'\big)},\\
    & u_V(\alpha',\kappa',n)=\frac{\Gamma\Big(\sum_{i=1}^d V\alpha_i'\Big)}{\Gamma\Big(n+V\sum_{i=1}^d \alpha_i'\Big)}\sum_{|a|=n}{n\choose a}\prod_{i=1}^d (\kappa_i')^{a_i}\frac{\Gamma(V\alpha_i'+a_i)}{\Gamma(V\alpha_i')}
\end{align*}
For $d=2$  we can also write the scaled $\tilde\pi_{n,V}$ using $a=(i,n-i)$ in terms of  Beta functions
\begin{equation}
     \tilde\pi_{n,V}(i)=\frac{1}{u_V(\alpha',\kappa',n)}(\kappa_1')^{i}(\kappa_2')^{(n-i)}{n \choose i} \frac{B\big(i+\frac{DV}{d\kappa_1'},n-i+\frac{DV}{d\kappa_2'}\big)}{B\big(\frac{DV}{d\kappa_1'},\frac{DV}{d\kappa_2'}\big)}
\end{equation}
as a weighted Beta-binomial distribution.
Just as in the symmetric case (see Sec 4.3 \cite{KBW20}), our proposed distribution $\tilde \Pi_V=\nu_V\tilde\pi_{n,V}$ will have a different shape: unimodal or multimodal, depending on whether: $DV<d$ or $DV>d$. 

\begin{figure}[hbt!]
\begin{subfigure}{}
  \centering
  \includegraphics[width=0.325\linewidth]{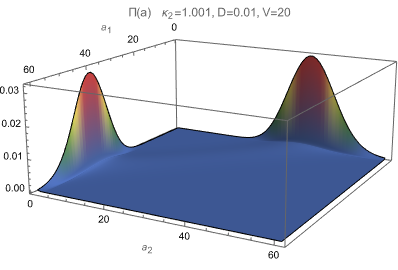}  
  \label{fig:sub-first}
\end{subfigure}
\begin{subfigure}{}
  \centering
  \includegraphics[width=0.325\linewidth]{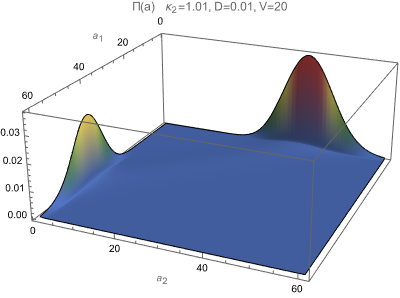}  
  \label{fig:sub-first}
\end{subfigure}
\begin{subfigure}{}
  \centering
  \includegraphics[width=0.325\linewidth]{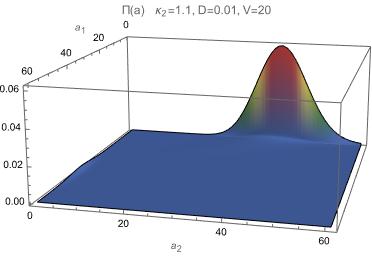}  
  \label{fig:sub-first}
\end{subfigure}

\caption{Plots of the analytic formula for $\tilde\Pi(a)$ when $DV=0.2<d=2$ show bimodal density. Parameters $D=0.01,V=20, \kappa_1'=1$  are fixed while $\kappa_2'$ is increased in plots from left to right: $\kappa_2'= 1.001, 1.01, 1.1$. Note the modes occur at $(a_1,a_2)\in\{(0,dV), (dV,0)\}$ with $dV=40$. }
\label{fig:dif k2}
\end{figure}

When $DV<d$ the parameters $\alpha'_1=DV/(d\kappa_1'), \alpha'_2=DV/(d\kappa_2')$ in the Beta-binomial are both less than 1, the weights have a small effect, and the distribution is bimodal. Since $\alpha'_2/\alpha_1'=\kappa_1'/\kappa_2'<1$, the mode at the value $(a_1,a_2)=(0,n)$ is higher than at $(a_1,a_2)=(n,0)$. Note that even very small asymmetries, e.g. $\kappa_2'=1.0001, \; \kappa'_1=1$, show a sizeable asymmetry in the size of the modes, and for a greater ratio, e.g. $\kappa_2'=1.1, \kappa'_1=1$ almost the entire mass of the distribution is concentrated at the higher mode (see Figure \ref{fig:dif k2}). Of course, on each plane $(0,n)$ and $(n,0)$ there is a mode  as well at $dV$ from the distribution of total molecular counts $\nu_V$.

\begin{figure}[hbt!]
\begin{subfigure}{}
  \centering
  \includegraphics[width=0.325\linewidth]{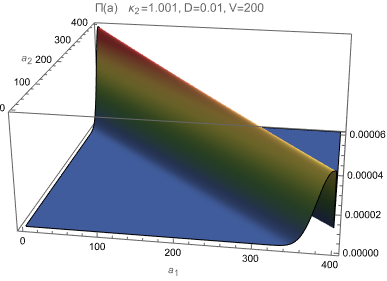}  
  \label{fig:sub-first}
\end{subfigure}
\begin{subfigure}{}
  \centering
  \includegraphics[width=0.325\linewidth]{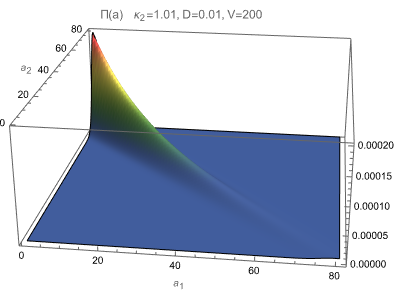}  
  \label{fig:sub-first}
\end{subfigure}
\begin{subfigure}{}
  \centering
  \includegraphics[width=0.325\linewidth]{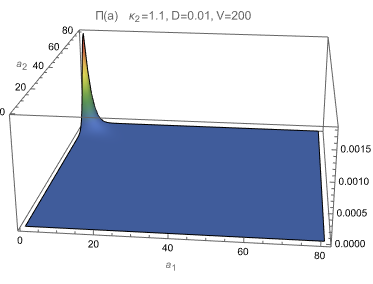}  
  \label{fig:sub-first}
\end{subfigure}

\caption{Plots of  $\tilde\Pi(a)$ when $DV=d=2$. Parameters are $D=0.01, V=200, \kappa_1'=1$, and $\kappa_2'$ increases in plots from left to right: $\kappa_2'=1.001, 1.01, 1.1$.}
\label{fig:propDV=2}
\end{figure}

In the intermediate case, when $DV=d=2$ we have $\alpha'_1=1$ and $\alpha_2'=1/\kappa_2'$. If $\kappa_2'=1$  the evaluation of $\tilde\Pi$ on  hyperplanes $E_n$ would be Uniform($[0,n]$), but for $\kappa_2'>1$ it  becomes skewed towards $\{(0,n)\}$, when e.g. $\kappa_2=1.1$ almost all the mass is at $(0,dV=400)$, see Figure \ref{fig:propDV=2}. 

\begin{figure}[hbt!]
\begin{subfigure}{}
  \centering
  \includegraphics[width=0.32\linewidth]{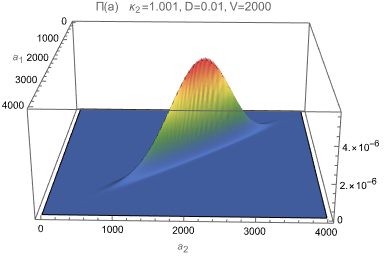}  
  \label{fig:sub-first}
\end{subfigure}
\begin{subfigure}{}
  \centering
  \includegraphics[width=0.32\linewidth]{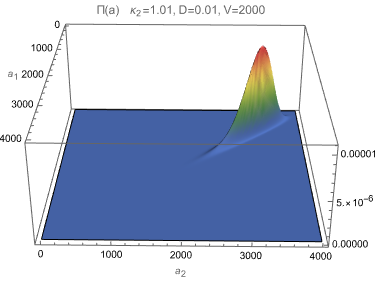}  
  \label{fig:sub-first}
\end{subfigure}
\begin{subfigure}{}
  \centering
  \includegraphics[width=0.32\linewidth]{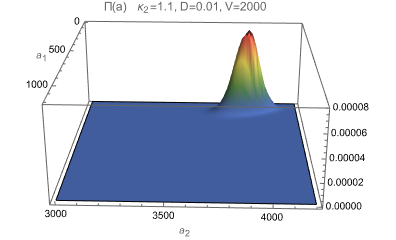}  
  \label{fig:sub-first}
\end{subfigure}
\caption{Plots of  $\tilde\Pi(a)$ when $DV=20>d=2$. Parameters are $D=0.01, V=2000, \kappa_1'=1$, and $\kappa_2'$ increases in plots from left to right: $\kappa_2'=1.001, 1.01, 1.1$.} 
\label{fig:propDV>2}
\end{figure}

When  $DV>2$ we have $\alpha'_1,\alpha'_2>1$ and $\tilde\Pi$ is unimodal concentrated around the fixed point $(a_1^*,a_2^*)$ of the deterministic system for this model
\begin{equation*}
\begin{split}
    \frac{d}{dt}{a_1}(t)&=(\kappa'_1-\kappa'_2)a_1a_2+\lambda'_1-\delta' a_1\\
    \frac{d}{dt}{a_2}(t)&=(\kappa'_2-\kappa'_1)a_1a_2+\lambda'_2-\delta' a_2
\end{split}
\end{equation*}
With $\kappa_1'=1$, increasing $\kappa_2'>1$ the mode moves such that the distribution is concentrated at a higher proportion of $a_2$ species to $a_1$ species. 
When $\kappa'_1=\kappa'_2$ this is a simple linear first order ODE with fixed point $(a_1^*,a_2^*)=(\lambda_1'/\delta',\lambda_2'/\delta')=(V,V).$
However, when $\kappa_1\neq\kappa_2$ this system of ODEs is nonlinear, and 
we found numerically for $\kappa_2'=1.001, 1.01, 1.1$ that the fixed points are, respectively, $(a_1^*,a_2^*)=(1801.96, 2198.04),  (763.93, 3236.07), (97.5, 3902.5)$ (Figure \ref{fig:propDV>2}, left to right).

\begin{figure}[hbt!]
\begin{subfigure}{}
  \centering
  \includegraphics[width=0.315\linewidth]{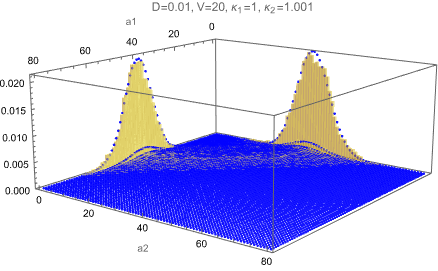}  
  \label{fig:sub-first}
\end{subfigure}
\begin{subfigure}{}
  \centering
  \includegraphics[width=0.315\linewidth]{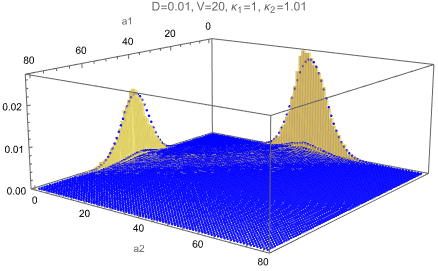}  
  \label{fig:sub-first}
\end{subfigure}
\begin{subfigure}{}
  \centering
  \includegraphics[width=0.315\linewidth]{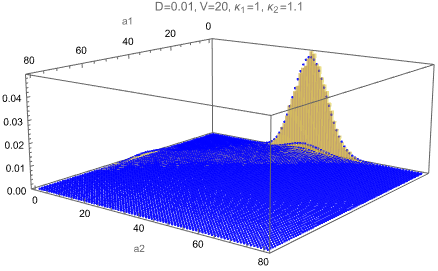}  
  \label{fig:sub-first}
\end{subfigure}
\caption{Histograms for  the fraction of occupation times in all states $a$ of $X$ are shown in yellow with $\tilde\Pi(a)$ overlaid in blue. The parameters are $ \lambda_1'=\lambda_2'=\delta'=D=0.01, V=20,  \kappa_1'=1$, varying $\kappa_2'=1.001,1.01,1.1$ from left to right.}
\label{fig:histplts}
\end{figure}

Since the Markov chain $(X_t)_{t\ge 0}$ is ergodic the long term fraction of its occupation times converges to its stationary distribution. We simulated the Markov chain, using the Gillespie algorithm, to further compare our approximate stationary distribution $\tilde\Pi$ with the occupation time when $DV<d$. 
In that case the modes at the boundaries of $E$ reflect the Discreteness Induced Transition (DIT) nature of the system, where most of the time is spent in the states $\{(0,n),(n,0)\}_{n\ge 0}$, and the quick transitions between these states results in the lack of occupation time spent in the interior of $E$.  As the asymmetry $\kappa_2'/\kappa_1'>1$ is increased, the DITs become more and more rare, as the autocatalytic reaction $A_1 +A_2 \xrightarrow{\kappa_2} 2A_2$ overpowers the autocatalytic reaction in the opposite direction.  Figure \ref{fig:histplts} shows the empirical long-term occupation time of the simulated Markov chain $X$ along with analytic plots of $\tilde\Pi$.  

\begin{figure}[hbt!]
\begin{subfigure}{}
  \centering
  \includegraphics[width=0.315\linewidth]{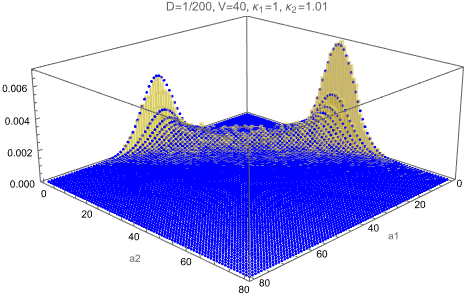}  
  \label{fig:sub-first}
\end{subfigure}
\begin{subfigure}{}
  \centering
  \includegraphics[width=0.325\linewidth]{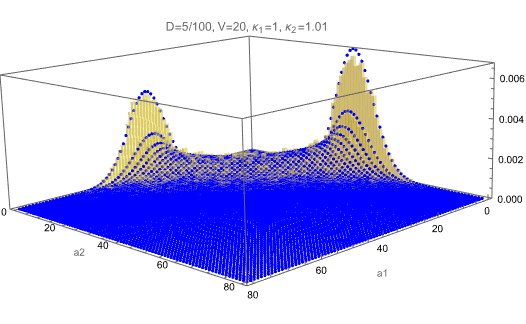}  
  \label{fig:sub-first}
\end{subfigure}
\begin{subfigure}{}
  \centering
  \includegraphics[width=0.3\linewidth]{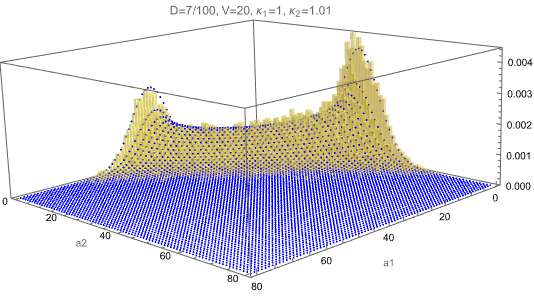}  
  \label{fig:sub-first}
\end{subfigure}

\caption{Histograms for  the fraction of occupation times in all states $a$ of $X$ are shown in yellow with $\tilde\Pi(a)$ overlaid in blue.  The parameters are $\lambda_1'=\lambda_2'=D, \kappa_1'=1, \kappa_2'=1.01 $, and $V=20$ was fixed while $D$ was varied  $D=2/100,5/100,7/100$ from left to right.}
\label{fig:histplst3}
\end{figure}

On the other hand, DITs still exist but are less pronounced if we keep $DV<d$  but we increase the rates of inflows $\lambda_i'$ and rate of outflows $\delta'$. In this case reactions $\emptyset  \xrightarrow{\lambda_i} A_i$, and $A_i \xrightarrow{\delta} \emptyset$ compete with the effects of auto-catalytic rates in either direction.  
The modes become less pronounced as the process spends an equal amount of time in the interior of $E$ as in states $\{(0,n),(n,0)\}_{n\ge 0}$. Figure \ref{fig:histplst3}  shows the effect of decreasing the inflow and outflow rates while keeping the volume fixed, and shows the empirical long term occupation time of the simulated Markov chain $X$ along with analytic plots of $\tilde\Pi$.  We note that despite the significant effect of the inflow and outflow reactions, our distribution $\tilde\Pi$ matches the empirical estimates of $\Pi$ well on all of $E$. This agrees with our numerical conclusions about $\mathcal{B}^*(a)\approx 0$ for all $a\in E$.

\section{Discussion}
In this paper, we studied a family of auto-catalytic reaction systems comprised of bimolecular interactions, plus inflows and outflows of molecular species from the system. The stochastic Markov chain model for the system is known to be ergodic with a unique stationary distribution, but the exact form of this distribution is only known in the case of symmetric auto-catalytic rates. To investigate the asymmetric case we proposed an approximate stationary distribution $\tilde\Pi$ in an explicit form of a weighted version of the distribution from the symmetric case. This distribution is known to be the stationary distribution of a related model in population genetics which is a conserved system with no inflows and outflows.

We derived an error function $\mathcal{B}^*$ which measures how close the proposed distribution is to the true stationary one. Numerical evaluations of this function, performed for a system with  $d=2$ species, show that our proposed distribution matches the true stationary distribution for almost all values in the state space (and is very close to it in the exceptional values). This is also confirmed by comparison of the empirical distribution of long-term fraction of occupation times for the Markov chain with analytic values of the proposed distribution.

The explicit form we provide for the proposed distribution $\tilde\Pi$ allows us to show that, under the usual mass-action scaling of the system using the order of magnitude of the overall species count as a parameter, the phenomenon of Discreteness Induced Transitions holds in the same parameter regime as in the symmetric case. The asymmetry of the auto-catalytic rates produces a significant asymmetry in the modes of the stationary distribution, even under very small rate differences. The magnitude of inflow and outflow rates relative to the auto-catalytic rates is shown to determine the extent to which Discreteness Induced Transitions are felt in the system, with the speed of transitions decreasing as the system becomes more open.

\subsubsection*{Acknowledgments}
This work was initiated following an American Institute of Mathematics (AIM) Workshop on Limits and Control of Stochastic Reaction Networks in August 2021, and partial results were presented at the AIM Square group meeting in October 2022. This work constitutes the MSc thesis project of C.G. ans was supported by L.P.'s NSERC Discovery Grant \# RGPIN-2015-06573.

\section{Appendices}\label{sec:appendix}
\subsection{Appendix A}
We start from 
\begin{align*}
&{\mathcal B}^*(a)=\frac{1}{\tilde\pi_n(a)}\big[ (L_{n-1}+L_n+L_{n+1})(a)-R_n(a)\big]\\
\mbox{with}&\\
    &R_n(a)=[\lambda_1+\lambda_2 +n\delta +(\kappa_1+\kappa_2)a_1a_2]\tilde\pi_n(a)\\
    &L_{n-1}(a)=\frac{n\delta \lambda_1}{\lambda_1+\lambda_2}\tilde\pi_{n-1}(a-e_1) + \frac{n\delta\lambda_2}{\lambda_1+\lambda_2}\tilde\pi_{n-1}(a-e_2)\\
    &L_n(a)=\kappa_2(a_1+1)(a_2-1)\tilde\pi_n(a+e_1-e_2)+\kappa_1(a_1-1)(a_2+1)\tilde\pi_n(a-e_1+e_2)\\
    &L_{n+1}(a)=\frac{\lambda_1+\lambda_2}{n+1}(a_1+1)\tilde\pi_{n+1}(a+e_1)+\frac{\lambda_1+\lambda_2}{n+1}(a_2+1)\tilde\pi_{n+1}(a+e_2).
\end{align*}
We use our expression \eqref{pid=2} (dropping the subscripts in the Gauss hypergeometric function $_2F_1$)
\begin{align*}
   \tilde \pi_n(i)&=\frac{1}{{ }F\big(-n,\alpha_1;1-\alpha_2-n;\frac{\kappa_1}{\kappa_2}\big)}\bigg(\frac{\kappa_1}{\kappa_2}\bigg)^a {n\choose i}  \frac{\Gamma(\alpha_1+i)\Gamma(\alpha_2+n-i)}{\Gamma(\alpha_1)\Gamma(\alpha_2+n)}
    \end{align*}
to get
\begin{align*}
    \frac{\tilde\pi_{n-1}(a-e_1)}{\tilde\pi_n(a)}&=\frac{a_1}{n}\frac{n-1+\alpha_2}{a_1-1+\alpha_1}\frac{\kappa_2}{\kappa_1}\frac{F(-n,\alpha_1,1-\alpha_2-n,\frac{\kappa_1}{\kappa_2})}{F(1-n,\alpha_1,2-\alpha_2-n,\frac{\kappa_1}{\kappa_2})}\\
    \frac{\tilde\pi_{n-1}(a-e_2)}{\tilde\pi_n(a)}&=\frac{a_2}{n}\frac{n-1+\alpha_2}{a_2-1+\alpha_2}\frac{F(-n,\alpha_1,1-\alpha_2-n,\frac{\kappa_1}{\kappa_2})}{F(1-n,\alpha_1,2-\alpha_2-n,\frac{\kappa_1}{\kappa_2})}\\
    \frac{\tilde\pi_n(a+e_1-e_2)}{\tilde\pi_n(a)}&=\frac{a_2}{a_1+1}\frac{a_1+\alpha_1}{a_2-1+\alpha_2}\frac{\kappa_1}{\kappa_2}\\
    \frac{\tilde\pi_n(a-e_1+e_2)}{\tilde\pi_n(a+e_1-e_2)}&=\frac{a_1}{a_2+1}\frac{a_2+\alpha_2}{a_1-1+\alpha_1}\frac{\kappa_2}{\kappa_1}\\
    \frac{\tilde\pi_{n+1}(a+e_1)}{\tilde\pi_n(a)}&=\frac{n+1}{a_1+1}\frac{a_1+\alpha_1}{n+\alpha_2}\frac{\kappa_1}{\kappa_2}\frac{F(-n,\alpha_1,1-\alpha_2-n,\frac{\kappa_1}{\kappa_2})}{F(-1-n,\alpha_1,-\alpha_2-n,\frac{\kappa_1}{\kappa_2})}\\
    \frac{\tilde\pi_{n+1}(a)}{\tilde\pi_n(a)}&=\frac{n+1}{a_2+1}\frac{a_2+\alpha_2}{n+\alpha_2}\frac{F(-n,\alpha_1,1-\alpha_2-n,\frac{\kappa_1}{\kappa_2})}{F(-1-n,\alpha_1,-\alpha_2-n,\frac{\kappa_1}{\kappa_2})}
\end{align*}
Then $L_{n-1},L_n,L_{n+1}$ divided by $\tilde\pi_n$ become
\begin{align*}
     &\frac{L_{n-1}}{\tilde\pi_n(a)}=\frac{(n-1+\alpha_2)\delta}{\lambda_1+\lambda_2}\frac{F(-n,\alpha_1,1-\alpha_2-n,\frac{\kappa_1}{\kappa_2})}{F(1-n,\alpha_1,2-\alpha_2-n,\frac{\kappa_1}{\kappa_2})}\Big[\frac{a_1\lambda_1}{a_1-1+\alpha_1}\frac{\kappa_2}{\kappa_1} +\frac{a_2\lambda_2}{a_2-1+\alpha_2}\Big]\\
     &\frac{L_n}{\tilde\pi_n(a)}=\kappa_1\frac{a_2-1}{a_2-1+\alpha_2}(a_1+\alpha_1)a_2 + \kappa_2\frac{a_1-1}{a_1-1+\alpha_1}(a_2+\alpha_2)a_1\\
     &\frac{L_{n+1}}{\tilde\pi_n(a)}=\frac{\lambda_1+\lambda_2}{(n+\alpha_2)}\frac{F(-n,\alpha_1,1-\alpha_2-n,\frac{\kappa_1}{\kappa_2})}{F(-1-n,\alpha_1,-\alpha_2-n,\frac{\kappa_1}{\kappa_2})}\Big[\frac{\kappa_1}{\kappa_2}(a_1+\alpha_1)+(a_2+\alpha_2)\Big]
\end{align*}
We use the fact that $\alpha_i=\frac{\delta}{\lambda_1+\lambda_2}\frac{\lambda_i}{\kappa_i}$ (so $\sum_{i=1}^2\kappa_i\alpha_i=\delta$), $a_1+a_2=n$, and express $\frac{L_{n-1}}{\tilde\pi_n}$ as
\begin{equation}\label{ln-1}
     \frac{L_{n-1}}{\tilde\pi_n(a)}=(n-1+\alpha_2)\frac{F(-n,\alpha_1,1-\alpha_2-n,\frac{\kappa_1}{\kappa_2})}{F(1-n,\alpha_1,2-\alpha_2-n,\frac{\kappa_1}{\kappa_2})}\kappa_2\Big[\frac{a_1\alpha_1}{a_1-1+\alpha_1}+\frac{a_2\alpha_2}{a_2-1+\alpha_2}\Big]
\end{equation}
For $\frac{L_n}{\tilde\pi_n(a)}$ we arrange terms as follows
\begin{equation*}
    \begin{split}
         \frac{L_n}{\tilde\pi_n(a)}=& \sum_{i=1}^2\sum_{j\neq i}\kappa_i(a_i+\alpha_i)\frac{a_j-1}{a_j-1+\alpha_j}a_j\\
         =& \sum_{i=1}^2\sum_{j\neq i}\kappa_i(a_i+\alpha_i)a_j-\sum_{i=1}^2\sum_{j\neq i}\frac{\kappa_i(a_i+\alpha_i)a_j\alpha_j}{a_j-1+\alpha_j}\\
          =& (\kappa_1+\kappa_2)a_1a_2 + \sum_{i=1}^2\sum_{j\neq i}\kappa_i\alpha_ia_j-\sum_{i=1}^2\sum_{j\neq i}\frac{\kappa_i(a_i+\alpha_i)a_j\alpha_j}{a_j-1+\alpha_j}\\
          =& (\kappa_1+\kappa_2)a_1a_2 + \sum_{i=1}^2\sum_{j=1}^2\kappa_i\alpha_ia_j-\sum_{i=1}^2\kappa_i\alpha_ia_i
          -\sum_{i=1}^2\sum_{j= 1}^2\frac{\kappa_j(a_j+\alpha_j)a_i\alpha_i}{a_i-1+\alpha_i}+\sum_{i= 1}^2\frac{\kappa_i(a_i+\alpha_i)a_i\alpha_i}{a_i-1+\alpha_i}\\
          =& (\kappa_1+\kappa_2)a_1a_2 + n\delta-\sum_{i=1}^2\kappa_i\alpha_ia_i
          -\sum_{i=1}^2\sum_{j= 1}^2\frac{\kappa_j(a_j+\alpha_j)a_i\alpha_i}{a_i-1+\alpha_i}+\sum_{i= 1}^2\kappa_ia_i\alpha_i+\sum_{i= 1}^2\frac{\kappa_ia_i\alpha_i}{a_i-1+\alpha_i}\\
          =& (\kappa_1+\kappa_2)a_1a_2 + n\delta-\sum_{i=1}^2\sum_{j= 1}^2\frac{\kappa_j(a_j+\alpha_j)a_i\alpha_i}{a_i-1+\alpha_i}+\sum_{i= 1}^2\frac{\kappa_ia_i\alpha_i}{a_i-1+\alpha_i}\\
          =& (\kappa_1+\kappa_2)a_1a_2 + n\delta-\Big(\sum_{j=1}^2\kappa_j(a_j+\alpha_j)\Big)\sum_{i=1}^2\frac{a_i\alpha_i}{a_i-1+\alpha_i} + \sum_{i= 1}^2\frac{\kappa_ia_i\alpha_i}{a_i-1+\alpha_i}\\
    \end{split}
\end{equation*}
so that the first two terms will cancel the last two terms in $\frac{R_n}{\tilde\pi_n}$; for the remaining terms we have
\begin{equation*}
    \begin{split}
        \quad &\sum_{i= 1}^2\frac{\kappa_ia_i\alpha_i}{a_i-1+\alpha_i} -\Big(\sum_{j=1}^2\kappa_j(a_j+\alpha_j)\Big)\sum_{i=1}^2\frac{a_i\alpha_i}{a_i-1+\alpha_i}\\
        &=\;\kappa_1 \frac{a_1\alpha_1}{a_1-1+\alpha_1}+\kappa_2\frac{a_2\alpha_2}{a_2-1+\alpha_2}\\
        &\;\;-\kappa_1(a_1+\alpha_1)\frac{a_1\alpha_1}{a_1-1+\alpha_1}-\kappa_2(a_2+\alpha_2)\frac{a_1\alpha_1}{a_1-1+\alpha_1}\\
        &\;\;-\kappa_1(a_1+\alpha_1)\frac{a_2\alpha_2}{a_2-2+\alpha_2}-\kappa_2(a_2+\alpha_2)\frac{a_2\alpha_2}{a_2-1+\alpha_2}\\
        &=\;\kappa_2\frac{a_1\alpha_1}{a_1-1+\alpha_1}\Big[\frac{\alpha_2}{\alpha_1}-\frac{\alpha_2}{\alpha_1}(a_1+\alpha_1)-(a_2+\alpha_2)\Big] 
        +\kappa_2\frac{a_2\alpha_2}{a_2-1+\alpha_2}\Big[1-\frac{\alpha_2}{\alpha_1}(a_1+\alpha_1)-(a_2+\alpha_2)\Big]\\
        &=\;\kappa_2\frac{a_1\alpha_1}{a_1-1+\alpha_1}\Big[\frac{\alpha_2}{\alpha_1}-\frac{\alpha_2}{\alpha_1}a_1-a_2-2\alpha_2\Big]+\kappa_2\frac{a_2\alpha_2}{a_2-1+\alpha_2}\Big[1-\frac{\alpha_2}{\alpha_1}a_1-a_2-2\alpha_2\Big]\\
        &=\;\kappa_2\frac{a_1\alpha_1}{a_1-1+\alpha_2}\Big[\frac{\alpha_2}{\alpha_1}+a_1-\frac{\alpha_2}{\alpha_1}a_1-1-\alpha_2\Big]+\kappa_2\frac{a_2\alpha_2}{a_2-1+\alpha_2}\Big[a_1-\frac{\alpha_2}{\alpha_1}a_1-\alpha_2\Big]\\
       &\;\;-\kappa_2(n-1+\alpha_2)\sum_{i=1}^2\frac{a_i\alpha_i}{a_i-1+\alpha_i}\\
       &=\;\kappa_2\sum_{i=1}^2\frac{a_i\alpha_i}{a_i-1+\alpha_i}\Big[a_1\big(1-\frac{\alpha_2}{\alpha_1}\big)-\alpha_2\Big]+\kappa_2\frac{a_1\alpha_1}{a_1-1+\alpha_1}\Big[\frac{\alpha_2}{\alpha_1}-1\Big]-\kappa_2(n-1+\alpha_2)\sum_{i=1}^2\frac{a_i\alpha_i}{a_i-1+\alpha_i}\\
        &=\;\sum_{i=1}^2\frac{a_i\alpha_i}{a_i-1+\alpha_i}\Big[a_1\big(\kappa_2-\kappa_1)-\kappa_2\alpha_2\Big]+\frac{a_1\alpha_1}{a_1-1+\alpha_1}\Big[\kappa_1-\kappa_2\Big]-\kappa_2(n-1+\alpha_2)\sum_{i=1}^2\frac{a_i\alpha_i}{a_i-1+\alpha_i}\\
    \end{split}
\end{equation*}
so that 
\begin{equation}\label{ln}
\begin{split}
    \frac{L_n}{\tilde\pi_n(a)}=&
    (\kappa_1+\kappa_2)a_1a_2+n\delta-\kappa_2(n-1+\alpha_2)\sum_{i=1}^2\frac{a_i\alpha_i}{a_i-1+\alpha_i}\\
    &+\sum_{i=1}^2\frac{a_i\alpha_i}{a_i-1+\alpha_i}\Big[a_1\big(\kappa_2-\kappa_1)-\alpha_2\Big]+\frac{a_1\alpha_1}{a_1-1+\alpha_1}\Big[\kappa_1-\kappa_2\Big]\\
    \end{split}
\end{equation}
where the first two terms will cancel with terms in $\frac{R_n}{\tilde\pi_n}$, and the third term is the same as in $\frac{L_{n-1}}{\tilde\pi_n}$ from \eqref{ln-1} apart from the ratio of Hypergeometric functions. \\
Finally, we express $\frac{L_{n+1}}{\tilde\pi_n}$ as
\begin{equation*}
    \begin{split}
        \frac{L_{n+1}}{\tilde\pi_n(a)}&=\frac{\lambda_1+\lambda_2}{n+\alpha_2}\frac{F(-n,\alpha_1,1-\alpha_2-n,\frac{\kappa_1}{\kappa_2})}{F(-1-n,\alpha_1,-\alpha_2-n,\frac{\kappa_1}{\kappa_2})}\Big[\frac{\kappa_1}{\kappa_2}(a_1+\alpha_1)+(a_2+\alpha_2)\Big]\\
        &=\frac{\lambda_1+\lambda_2}{n+\alpha_2}\frac{F(-n,\alpha_1,1-\alpha_2-n,\frac{\kappa_1}{\kappa_2})}{F(-1-n,\alpha_1,-\alpha_2-n,\frac{\kappa_1}{\kappa_2})}\Big[(n+\alpha_2)+\frac{\kappa_1}{\kappa_2}(a_1+\alpha_1)-a_1\Big]\\ 
        &=(\lambda_1+\lambda_2)\frac{F(-n,\alpha_1,1-\alpha_2-n,\frac{\kappa_1}{\kappa_2})}{F(-1-n,\alpha_1,-\alpha_2-n,\frac{\kappa_1}{\kappa_2})}
        +\frac{\lambda_1+\lambda_2}{n+\alpha_2}\frac{F(-n,\alpha_1,1-\alpha_2-n,\frac{\kappa_1}{\kappa_2})}{F(-1-n,\alpha_1,-\alpha_2-n,\frac{\kappa_1}{\kappa_2})}\Big[a_1\big(\frac{\kappa_1}{\kappa_2}-1\big)+\alpha_2\Big]
        \end{split}
\end{equation*}
For ${\mathcal B}^*=\frac{1}{\tilde\pi_n}(L_{n-1}+L_n-L_{n-1}-R_n)$ we now get the expression in \eqref{DBeqn}
  \begin{equation*}
 \begin{split}
     {\mathcal B}^*(a)=&(\lambda_1+\lambda_2)\Bigg[1-\frac{F(-n,\alpha_1,1-\alpha_2-n,\frac{\kappa_1}{\kappa_2})}{F(-1-n,\alpha_1,-\alpha_2-n,\frac{\kappa_1}{\kappa_2})}\Bigg]\\
     &-\kappa_2(n-1+\alpha_2)\sum_{i=1}^2\frac{a_i\alpha_i}{a_i-1+\alpha_i}\Bigg[\frac{F(-n,\alpha_1,1-\alpha_2-n,\frac{\kappa_1}{\kappa_2})}{F(1-n,\alpha_1,2-\alpha_2-n,\frac{\kappa_1}{\kappa_2})}-1\Bigg]\\
      &-\sum_{i=1}^2\frac{a_i\alpha_i}{a_i-1+\alpha_i}\Big[a_1\big(\kappa_2-\kappa_1)-\kappa_2\alpha_2\Big]-\frac{a_1\alpha_1}{a_1-1+\alpha_1}\Big[\kappa_1-\kappa_2\Big]\\
      &-\frac{\lambda_1+\lambda_2}{n+\alpha_2}\frac{F(-n,\alpha_1,1-\alpha_2-n,\frac{\kappa_1}{\kappa_2})}{F(-1-n,\alpha_1,-\alpha_2-n,\frac{\kappa_1}{\kappa_2})}\Big[a_1\big(\frac{\kappa_1}{\kappa_2}-1\big)+\alpha_2\Big].   
\end{split}
 \end{equation*}

\subsection{Appendix B}
We strat from \eqref{DBeqn} and use the \eqref{scalingKBW} scaling by $V$ 
\begin{equation*}
\kappa_i=\frac{\kappa_i'}{V},\;\; \lambda_i=DV,\;\; \delta=D,\;\;  \quad\Rightarrow \quad \alpha_i=\frac{DV}{d\kappa_i'}=V\alpha_i'
\end{equation*}
 to get 
\begin{equation*}
 \begin{split}
     {\mathcal B}^*_V(a)=
     &dDV\Bigg[1-\frac{F(-n,V\alpha_1',1-V\alpha_2'-n,\frac{\frac{\kappa_1'}{V}}{\frac{\kappa_2'}{V}})}{F(-1-n,V\alpha_1',-V\alpha_2'-n,\frac{\frac{\kappa_1'}{V}}{\frac{\kappa_2'}{V}})}\Bigg]\\
     &-\frac{\kappa_2'}{V}(n-1+V\alpha_2')\sum_{i=1}^2\frac{a_iV\alpha_i'}{a_i-1+V\alpha_i'}\Bigg[\frac{F(-n,V\alpha_1',1-V\alpha_2'-n,\frac{\frac{\kappa_1'}{V}}{\frac{\kappa_2'}{V}})}{F(1-n,V\alpha_1',2-V\alpha_2'-n,\frac{\frac{\kappa_1'}{V}}{\frac{\kappa_2'}{V}})}-1\Bigg]\\
      &-\sum_{i=1}^2\frac{a_iV\alpha_i'}{a_i-1+VV\alpha_i''}\Big[a_1\big(\frac{\kappa_2'}{V}-\frac{\kappa_1'}{V})-\frac{\kappa_2'}{V}V\alpha_2'\Big]-\frac{a_1V\alpha_1'}{a_1-1+V\alpha_1'}\Big[\frac{\kappa_1'}{V}-\frac{\kappa_2'}{V}\Big]\\
      &-\frac{\lambda_1+\lambda_2}{n+V\alpha_2'}\frac{F(-n,V\alpha_1',1-V\alpha_2'-n,\frac{\frac{\kappa_1'}{V}}{\frac{\kappa_2'}{V}})}{F(-1-n,V\alpha_1',-V\alpha_2'-n,\frac{\frac{\kappa_1'}{V}}{\frac{\kappa_2'}{V}})}\Big[a_1\big(\frac{\frac{\kappa_1'}{V}}{\frac{\kappa_2'}{V}}-1\big)+V\alpha_2'\Big]\\   
\end{split}
 \end{equation*}
 \begin{equation*}
 \begin{split}
     =&dDV\Bigg[1-\frac{F(-n,V\alpha_1',1-V\alpha_2'-n,\frac{\kappa_1'}{\kappa_2'})}{F(-1-n,V\alpha_1',-V\alpha_2'-n,\frac{\kappa_1'}{\kappa_2'})}\Bigg]\\
     &-\kappa_2(n-1+V\alpha_2')\sum_{i=1}^2\frac{a_i\alpha_i'}{a_i-1+V\alpha_i'}\Bigg[\frac{F(-n,V\alpha_1',1-V\alpha_2'-n,\frac{\kappa_1'}{\kappa_2'})}{F(1-n,V\alpha_1',2-V\alpha_2'-n,\frac{\kappa_1'}{\kappa_2'})}-1\Bigg]\\
      &-\sum_{i=1}^2\frac{a_i\alpha_i'}{a_i-1+V\alpha_i'}\Big[a_1\big(\kappa_2-\kappa_1)-\kappa_2\alpha_2'\Big]-\frac{a_1\alpha_1'}{a_1-1+V\alpha_1'}\Big[\kappa_1'-\kappa_2'\Big]\\
      &-\frac{dDV}{n+V\alpha_2'}\frac{F(-n,V\alpha_1',1-V\alpha_2'-n,\frac{\kappa_1'}{\kappa_2'})}{F(-1-n,V\alpha_1',-V\alpha_2'-n,\frac{\kappa_1'}{\kappa_2'})}\Big[a_1\big(\frac{\kappa_1'}{\kappa_2'}-1\big)+V\alpha_2'\Big]\\   
\end{split}
 \end{equation*}
Substituting  $\alpha_i'=\frac{D}{d\kappa_i'}$ and rearranging terms we get the expression in \eqref{BV}
 \begin{equation*}
 \begin{split}
      {\mathcal B}^*_V(a) =&dDV\Bigg(\Bigg[1-\frac{F(-n,\frac{DV}{d\kappa_1'},1-\frac{DV}{d\kappa_2'}-n,\frac{\kappa_1'}{\kappa_2'})}{F(-1-n,\frac{DV}{d\kappa_1'},-\frac{DV}{d\kappa_2'}-n,\frac{\kappa_1'}{\kappa_2'})}\Bigg]\\
     +&\frac{1}{n+\frac{DV}{d\kappa_2'}}\frac{F(-n,\frac{DV}{d\kappa_1'},1-\frac{DV}{d\kappa_2'}-n,\frac{\kappa_1'}{\kappa_2'})}{F(-1-n,\frac{DV}{d\kappa_1'},-\frac{DV}{d\kappa_2'}-n,\frac{\kappa_1'}{\kappa_2'})}\Big[a_1\big(1-\frac{\kappa_1'}{\kappa_2'}\big)-\frac{DV}{d\kappa_2'}\Big]\Bigg)\\
     -&\kappa_2'(n-1+\frac{DV}{d\kappa_2'})\Bigg[\frac{F(-n,\frac{DV}{d\kappa_1'},1-\frac{DV}{d\kappa_2'}-n,\frac{\kappa_1'}{\kappa_2'})}{F(1-n,\frac{DV}{d\kappa_1'},2-\frac{DV}{d\kappa_2'}-n,\frac{\kappa_1'}{\kappa_2'})}-1\Bigg]\bigg(\frac{a_1\frac{DV}{d\kappa_1'}}{a_1-1+\frac{DV}{d\kappa_1'}}+\frac{a_2\frac{DV}{d\kappa_2'}}{a_2-1+\frac{DV}{d\kappa_2'}}\bigg)\\
      -&\kappa_2'\Big[a_1\big(1-\frac{\kappa_1'}{\kappa_2'})-\frac{DV}{d}\Big]\bigg(\frac{a_1\frac{DV}{d\kappa_1'}}{a_1-1+\frac{DV}{d\kappa_1'}}+\frac{a_2\frac{DV}{d\kappa_2'}}{a_2-1+\frac{DV}{d\kappa_2'}}\bigg)+\kappa_2'\frac{a_1\frac{D}{d\kappa_1'}}{a_1-1+\frac{DV}{d\kappa_1'}}\Big[1-\frac{\kappa_1'}{\kappa_2'}\Big].
\end{split}
 \end{equation*}



\end{document}